\documentclass[10pt]{article}

\usepackage{amsmath}
\usepackage{amsfonts}
\usepackage{amssymb}
\usepackage[english]{babel}
\usepackage{amsthm}

\def\pf{\par\noindent {\bf Proof}~\par\noindent}
\def\qed{~\hfill{$\square$}\pagebreak[1]\par\medskip\par}

\newcommand{\mR}{\mathbb{R}}
\newcommand{\mC}{\mathbb{C}}

\newcommand{\mF}{\mathbb{F}}
\newcommand{\p}{\partial}
\newcommand{\ux}{\underline{x}}

\newcommand{\uom}{\underline{\omega}}
\newcommand{\pux}{\underline{\partial}}
\newcommand{\po}{\partial_{x_0}}

\newcommand{\mcH}{\mathcal{H}}

\newcommand{\enb}{\overline{e_0}}
\newcommand{\Dbar}{\overline{D}}
\newcommand{\pxn}{\partial_{x_0}}
\newcommand{\onehalf}{\frac{1}{2}}

\newcommand{\dkw}{(x_0^2 + r^2)}
\newcommand{\dk}{x_0^2 + r^2}
\newcommand{\dkwo}{(x_0^2 + \rho^2)}
\newcommand{\dko}{x_0^2 + \rho^2}

\newtheorem{proposition}{Proposition}
\newtheorem{remark}{Remark}

\newtheorem{property}{Property}

\numberwithin{equation}{section}
\numberwithin{theorem}{section}
\numberwithin{proposition}{section}
\numberwithin{lemma}{section}
\numberwithin{definition}{section}
\numberwithin{remark}{section}
\numberwithin{corollary}{section}

\begin{document}

\title{Harmonic and Monogenic Potentials in Low Dimensional Euclidean Half--Space}

\author{F.\ Brackx, H.\ De Bie, H.\ De Schepper}

\date{\small{Clifford Research Group, Department of Mathematical Analysis,\\ Faculty of Engineering and Architecture, Ghent University\\
Building S22, Galglaan 2, B-9000 Gent, Belgium\\}}

\maketitle

\begin{abstract}
\noindent In the framework of Clifford analysis, a chain of harmonic and monogenic potentials in the upper half of Euclidean space $\mR^{m+1}$ was constructed recently, including a higher dimensional analogue of the logarithmic function in the complex plane. Their distributional limits at the boundary $\mR^{m}$ were also determined. In this paper the potentials and their distributional boundary values are calculated in dimensions 3 and 4, dimensions for which the expressions in general dimension break down.
\end{abstract}

\maketitle


\section{Introduction}
\label{intro}


Recently, see \cite{bdbds1, bdbds2}, a generalization to Euclidean upper half--space $\mR^{m+1}_+$ was constructed of the  logarithmic function $\ln{z}$ which is holomorphic in the upper half of the complex plane. This construction was carried out in the framework of Clifford analysis, where the functions under consideration take their values in the universal Clifford algebra $\mR_{0,m+1}$ constructed over the Euclidean space $\mR^{m+1}$ equipped with a quadratic form of signature $(0,m+1)$. The concept of a higher dimensional holomorphic function, mostly called {\em monogenic} function, is expressed by means of the generalized Cauchy--Riemann operator $D$, which is a combination of the derivative with respect to one of the real variables, say $x_0$, and the so--called Dirac operator $\pux$ in the remaining real variables $(x_1, x_2, \ldots, x_m)$. This generalized Cauchy--Riemann operator $D$ and its Clifford algebra conjugate $\overline{D}$ linearize the Laplace operator, whence Clifford analysis may be seen as a refinement of harmonic analysis.\\[-2mm]

The starting point for constructing a higher dimensional monogenic logarithmic function, was the fundamental solution of the generalized Cauchy--Riemann operator $D$, also called Cauchy kernel, and its relation to the Poisson kernel and its harmonic conjugate in $\mR^{m+1}_+$. We then proceeded by induction in two directions, {\em downstream} by differentiation and {\em upstream} by primitivation, yielding a doubly infinite chain of monogenic, and thus harmonic, potentials. This chain mimics the well--known sequence of holomorphic potentials in $\mC_+$ (see e.g. \cite{slang}):
$$
\frac{1}{k!} z^k \left[ \ln z - ( 1 + \frac{1}{2} + \ldots + \frac{1}{k}) \right] \rightarrow \ldots \rightarrow z ( \ln z - 1) \rightarrow \ln z 
\stackrel{\frac{d}{dz}}{\longrightarrow} \frac{1}{z} \rightarrow - \frac{1}{z^2} \rightarrow \ldots \rightarrow (-1)^{k-1} \frac{(k-1)!}{z^k}
$$
Identifying the boundary of upper half--space with $\mR^m \cong \{(x_0,\ux) \in \mR^{m+1} : x_0 = 0\}$, the distributional limits for $x_0 \rightarrow 0+$ of those potentials were computed. They split up into two classes of distributions, which are linked by the Hilbert transform, one scalar--valued, the second one Clifford vector--valued. They belong to two out of the four families of Clifford distributions which were thoroughly studied in a series of papers, see \cite{fb1, fb2, distrib} and the references therein.\\

The general expressions for the potentials and their boundary values, as established in \cite{bdbds1, bdbds2}, break down for dimensions $3 \ (m=2)$ and $4 \ (m=3)$, which are particularly important with a view on applications. In these specific cases ad hoc calculations have to be carried out. This is the aim of the underlying paper, which is a valuable and useful complement to  \cite{bdbds1, bdbds2} . To make the paper self--contained the basics of Clifford algebra and Clifford analysis are recalled.


\section{Basics of Clifford analysis}
\label{basics}


Clifford analysis (see e.g. \cite{red, green}) is a function theory which offers a natural and elegant generalization to higher dimension of holomorphic functions in the complex plane and refines harmonic analysis. Let  $(e_0, e_1,\ldots,e_m)$ be the canonical orthonormal basis of Euclidean space $\mR^{m+1}$ equipped with a quadratic form of signature $(0,m+1)$. Then the non--commutative multiplication in the universal real Clifford algebra $\mR_{0,m+1}$ is governed by the rule 
$$
e_{\alpha} e_{\beta} + e_{\beta} e_{\alpha} = -2 \delta_{\alpha \beta}, \qquad \alpha,\beta = 0, 1,\ldots,m
$$
whence $\mR_{0,m+1}$ is generated additively by the elements $e_A = e_{j_1} \ldots e_{j_h}$, where $A=\lbrace j_1,\ldots,j_h \rbrace \subset \lbrace 0,\ldots,m \rbrace$, with $0\leq j_1<j_2<\cdots < j_h \leq m$, and $e_{\emptyset}=1$. 
For an account on Clifford algebra we refer to e.g. \cite{porteous}.\\[-2mm]

We identify the point $(x_0, x_1, \ldots, x_m) \in \mR^{m+1}$ with the Clifford--vector variable 
$$
x = x_0 e_0  + x_1 e_1  + \cdots x_m e_m = x_0 e_0  + \ux
$$ 
and the point $(x_1, \ldots, x_m) \in \mR^{m}$ with the Clifford--vector variable $\ux$. 
The introduction of spherical co--ordinates $\ux = r \uom$, $r = |\ux|$, $\uom \in S^{m-1}$, gives rise to the Clifford--vector valued locally integrable function $\uom$, which is to be seen as the higher dimensional analogue of the {\em signum}--distribution on the real line.\\[-2mm]

At the heart of Clifford analysis lies the so--called Dirac operator 
$$
\p = \pxn e_0 + \p_{x_1} e_1 + \cdots \p_{x_m} e_m =  \pxn e_0 + \pux
$$
which squares to the negative Laplace operator: $\p^2 = - \Delta_{m+1}$, while also $\pux^2 = - \Delta_{m}$. The fundamental solution of the Dirac operator $\p$ is given by (see \cite{red,green})
$$
E_{m+1} (x) = - \frac{1}{\sigma_{m+1}} \ \frac{x}{|x|^{m+1}}
$$
where $\sigma_{m+1} = \frac{2\pi^{\frac{m+1}{2}}}{\Gamma(\frac{m+1}{2})}$ stands for the area of the unit sphere $S^{m}$ in $\mR^{m+1}$. It thus holds
$$
\p E_{m+1} (x) = \delta(x)
$$
with $\delta(x)$ the standard Dirac distribution in $\mR^{m+1}$.\\
We also introduce the generalized Cauchy--Riemann operator 
$$
D = \onehalf \enb \p = \onehalf (\pxn + \enb \pux)
$$ 
which, together with its Clifford algebra conjugate $\Dbar = \onehalf(\pxn - \enb \pux)$, linearizes the Laplace operator: $D \Dbar = \Dbar D = \frac{1}{4} \Delta_{m+1}$. \\[-2mm]

A continuously differentiable function $F(x)$, defined in an open region $\Omega \subset \mR^{m+1}$ and taking values in the Clifford algebra  $\mR_{0,m+1}$ (or subspaces thereof), is called (left--)monogenic if it satisfies  in $\Omega$ the equation $D F = 0$, which is equivalent with $\p F = 0$. \\

Singling out the basis vector $e_0$, we can decompose the real Clifford algebra $\mR_{0,m+1}$ in terms of the Clifford algebra 
$\mR_{0,m}$ as $\mR_{0,m+1} = \mR_{0,m} \oplus \enb \, \mR_{0,m}$. Similarly we decompose the considered functions as 
$$
F(x_0,\ux) = F_1(x_0,\ux) + \enb \, F_2(x_0,\ux)
$$ 
where $F_1$ and $F_2$ take their values in the Clifford algebra $\mR_{0,m}$; mimicking functions of a complex variable, we will call $F_1$ the {\em real} part and $F_2$ the {\em imaginary} part of the function $F$.\\[-2mm]


\section{Harmonic and monogenic potentials in $\mR_+^{3}$}
\label{potentials3}


In this section we calculate the harmonic and monogenic potentials in upper half--space $\mR^{3}_+$, as defined in general in \cite{bdbds1, bdbds2}.

\subsection{The Cauchy kernel in $\mR_+^{3}$}

The starting point is the Cauchy kernel of Clifford analysis, i.e. the fundamental solution of the generalized Cauchy--Riemann operator $D$, which is, up to a constant factor, the fundamental solution $E_3(x)$ of the Dirac operator $\p$, and is  given by
$$
C_{-1}(x_0,\ux) = \frac{1}{\sigma_{3}} \, \frac{x \overline{e_0}}{|x|^{3}} =  \frac{1}{4\pi} \, \frac{x_0 - \overline{e_0} \ux}{(x_0^2+r^2)^{3/2}}
$$
where we have put $|\ux| = r$.
It may be decomposed in terms of the traditional Poisson kernel $P(x_0,\ux)$ and its conjugate $Q(x_0,\ux)$ in $\mR^{3}_+$:
$$
C_{-1}(x_0,\ux) = \frac{1}{2} A_{-1}(x_0,\ux) + \frac{1}{2} \overline{e_0} \, B_{-1}(x_0,\ux)
$$
where for $x_0 >0$,
\begin{eqnarray*}
A_{-1}(x_0,\ux) & = & P(x_0,\ux) \ = \ \phantom{-} \frac{1}{2\pi} \, \frac{x_0}{(x_0^2+r^2)^{3/2}} \label{A-1}\\
B_{-1}(x_0,\ux) & = & Q(x_0,\ux) \ = \ - \frac{1}{2\pi} \, \frac{\ux}{(x_0^2+r^2)^{3/2}} \label{B-1}
\end{eqnarray*}
Their distributional limits for $x_0 \rightarrow 0+$ are given by
\begin{eqnarray*}
a_{-1}(\ux) & = & \lim_{x_0 \rightarrow 0+} A_{-1}(x_0,\ux) \ = \ \delta(\ux)\\
b_{-1}(\ux) & = & \lim_{x_0 \rightarrow 0+} B_{-1}(x_0,\ux) \ = \ H(\ux) \ = \ - \frac{1}{2\pi} \, \mbox{Pv} \frac{\ux}{r^3} \ = \ - \frac{1}{2\pi} \, \mbox{Pv} \frac{\uom}{r^2}
\end{eqnarray*}
with $\mbox{Pv}$ standing for the "principal value" distribution in $\mR^2$. The distribution
$
H(\ux)
$
is the convolution kernel of the Hilbert transform $\mathcal{H}$ in $\mR^2$ (see e.g.\ \cite{gilmur}). Note that $\mathcal{H}^2 = \mathbf{1}$ and that this Hilbert transform links both distributional boundary values, as expressed in the following property, which in fact holds regardless the dimension.
\begin{property} One has
$$
\begin{array}{rc}
(i)  \quad \mathcal{H} \left [ a_{-1} \right ]  =  \mathcal{H} \left [ \delta \right ] \ = \ H \ast \delta \ = \ H \ = \ b_{-1} \\[2mm]
(ii) \quad  \mathcal{H} \left [ b_{-1} \right ]  =  \mathcal{H} \left [ H \right ] \ = \ H \ast H \ = \ \delta \ = \ a_{-1}
\end{array}
$$
\end{property}
\noindent Note also that $a_{-1}(\ux) = \delta(\ux) = E_0$ can be seen as the fundamental solution of the identity operator $\pux^0 = {\bf 1}$, while $b_{-1}(\ux) = H(\ux) =  - \frac{1}{2\pi} \, \mbox{Pv} \frac{\uom}{r^2} = F_0$ is the fundamental solution of the Hilbert operator $^0\mcH = \mcH$ (see \cite{bdbds2}).

\subsection{The downstream potentials in $\mR_+^{3}$}

The first in the sequence of so--called {\em downstream} potentials is the function $C_{-2}$ defined by
$$
\overline{D} C_{-1} = C_{-2} = \frac{1}{2} A_{-2} + \frac{1}{2} \overline{e_0} B_{-2}
$$
Clearly it is monogenic  in $\mR^{3}_+$, since $DC_{-2} = D \overline{D} C_{-1} = \frac{1}{4} \Delta_{m+1} C_{-1} = 0$. 
The definition itself of $C_{-2}(x_0,\ux)$ implies that it shows the monogenic potential (or primitive) $C_{-1}(x_0,\ux)$ and the conjugate harmonic potentials $A_{-1}(x_0,\ux)$ and $\overline{e_0} B_{-1}(x_0,\ux)$. For the notion of higher dimensional {\em conjugate harmonicity} in the framework of Clifford analysis we refer the reader to \cite{fbrdfs}.\\

The harmonic component $A_{-2}(x_0,\ux)$ may be calculated as $\po A_{-1}$ or as $-\pux B_{-1}$, and even the general expression is valid, leading to
$$
A_{-2}(x_0,\ux) = \frac{1}{2\pi} \, \frac{-2x_0^2 + r^2}{\dkw^{5/2}}
$$
Similarly the harmonic component $B_{-2}(x_0,\ux)$ may be calculated as $\po B_{-1}$ or as $-\pux A_{-1}$, and also the general expression is valid, leading to
$$
B_{-2}(x_0,\ux) = \frac{1}{2\pi} \, \frac{3 x_0 \ux}{\dkw^{5/2}}
$$
The distributional limits for $x_0 \rightarrow 0+$ of these harmonic potentials are given by
$$
\left \{ \begin{array}{rcl}
a_{-2}(\ux) = \lim_{x_0 \rightarrow 0+} A_{-2}(x_0,\ux) & = & \displaystyle\frac{1}{2\pi} \, {\rm Fp} \displaystyle\frac{1}{r^3}\\[4mm]
b_{-2}(\ux) = \lim_{x_0 \rightarrow 0+} B_{-2}(x_0,\ux) & = & - \pux \delta
\end{array} \right .
$$
where ${\rm Fp}$ stands for the "finite part" distribution on the real $r$--line.\\
Note that $a_{-2}(\ux) = - F_{-1}$, with $F_{-1} = \pux H$ the fundamental solution of the operator $^{-1}\mcH$, while 
$b_{-2}(\ux) = - E_{-1}$, with $E_{-1} = \pux \delta$ the fundamental solution of the operator $\pux^{-1}$ (see \cite{bdbds2}).\\

Proceeding in the same manner, the sequence of {\em downstream} monogenic potentials in $\mR_+^{m+1}$ is defined by
$$
C_{-k-1} = \overline{D} C_{-k} = \overline{D}^2 C_{-k+1} = \ldots = \overline{D}^k C_{-1}, \qquad k=1,2,\ldots
$$
where each monogenic potential decomposes into two conjugate harmonic potentials:
$$
C_{-k-1} = \frac{1}{2} A_{-k-1} + \frac{1}{2} \overline{e_0} B_{-k-1}, \qquad k=1,2,\ldots
$$
with, for $k$ odd, say $k=2\ell-1$,
$$
\left \{ \begin{array}{rcl}
A_{-2\ell} & = & \p_{x_0}^{2\ell-1} A_{-1} \ = \ - \p_{x_0}^{2\ell-2} \pux B_{-1} \ = \ \ldots \ = \ - \pux^{2\ell-1} B_{-1} \\[2mm]
B_{-2\ell} & = & \p_{x_0}^{2\ell-1} B_{-1} \ = \ - \p_{x_0}^{2\ell-2} \pux A_{-1} \ = \ \ldots \ = \ - \pux^{2\ell-1} A_{-1} 
\end{array} \right .
$$
while for $k$ even, say $k=2\ell$,
$$
\left \{ \begin{array}{rcl}
A_{-2\ell-1} & = & \p_{x_0}^{2\ell} A_{-1} \ = \ - \p_{x_0}^{2\ell-1} \pux B_{-1} \ = \ \ldots \ = \ \pux^{2\ell} A_{-1} \\[2mm]
B_{-2\ell-1} & = & \p_{x_0}^{2\ell} B_{-1} \ = \ - \p_{x_0}^{2\ell-1} \pux A_{-1} \ = \ \ldots \ = \ \pux^{2\ell} B_{-1} 
\end{array} \right .
$$
The harmonic potentials $A_{-k}(x_0,\ux)$ are real--valued and given by
$$
A_{-k}(x_0,\ux) = (-1)^{k+1} \, \frac{1}{2\pi} \, k! \, \frac{1\cdot3\cdot\cdots(2k-1)}{(k+1)(k+2)\cdots(2k)} \, \frac{r^k}{\dkw^{k+\onehalf}} \, i^k \, C_k^{-k}(i \frac{x_0}{r})
$$
where $C_k^{\lambda}$ stands for the Gegenbauer polynomial. We have explicitly calculated $A_{-k}$ for the $k$--values $1,2,3$ and $4$. We obtain\\

\noindent
for $k=1$
$$
A_{-1} = \frac{1}{2\pi} \, \onehalf \, \frac{r}{\dkw^{3/2}} \, i \, C_1^{-1}(i \frac{x_0}{r}) = \frac{1}{2\pi} \,  \frac{x_0}{\dkw^{3/2}}
$$
for $k=2$
$$
A_{-2} = - \frac{1}{2\pi} \, 2! \, \frac{1\cdot3}{3\cdot4} \, \frac{r^2}{\dkw^{5/2}} \, i^2 \, C_2^{-2}(i \frac{x_0}{r}) = \frac{1}{2\pi} \,  \frac{-2 x_0^2 + r^2}{\dkw^{5/2}}
$$
for $k=3$
$$
A_{-3} =  \frac{1}{2\pi} \, 3! \, \frac{1\cdot3\cdot5}{4\cdot5\cdot6} \, \frac{r^3}{\dkw^{7/2}} \, i^3 \, C_3^{-3}(i \frac{x_0}{r}) = \frac{1}{2\pi} \,  \frac{6 x_0^3 - 9 x_0 r^2}{\dkw^{7/2}}
$$
and for $k=4$
$$
A_{-4} = - \frac{1}{2\pi} \, 4! \, \frac{1\cdot3\cdot5\cdot7}{5\cdot6\cdot7\cdot8} \, \frac{r^4}{\dkw^{9/2}} \, i^4 \, C_4^{-4}(i \frac{x_0}{r}) = \frac{1}{2\pi} \,  \frac{-24 x_0^4 + 72  x_0^2  r^2 - 9 r^4}{\dkw^{9/2}}
$$
The harmonic potentials $B_{-k}(x_0,\ux)$ are Clifford vector--valued and given by
$$
B_{-k}(x_0,\ux) = (-1)^{k} \, \frac{1}{2\pi} \, (k-1)! \, \frac{3\cdot5\cdots(2k-1)}{(k+2)(k+3)\cdots(2k)} \, \frac{r^{k-1} \ux}{\dkw^{k+\onehalf}} \, i^{k-1} \, C_{k-1}^{-k}(i \frac{x_0}{r})
$$
We have explicitly calculated $B_{-k}$ for the $k$--values $1,2,3$ and $4$. We obtain\\

\noindent
for $k=1$
$$
B_{-1} = - \frac{1}{2\pi} \, \frac{\ux}{\dkw^{3/2}} \, C_0^{-1}(i \frac{x_0}{r}) = - \frac{1}{2\pi} \,  \frac{\ux}{\dkw^{3/2}}
$$
for $k=2$
$$
B_{-2} =  \frac{1}{2\pi} \, \frac{3}{4} \, \frac{r \ux}{\dkw^{5/2}} \, i \, C_1^{-2}(i \frac{x_0}{r}) = \frac{1}{2\pi} \,  \frac{3 x_0 \ux}{\dkw^{5/2}}
$$
for $k=3$
$$
B_{-3} =  - \frac{1}{2\pi} \, 2! \,  \frac{3\cdot5}{5\cdot6} \, \frac{r^2 \ux}{\dkw^{7/2}} \, i^2 \, C_2^{-3}(i \frac{x_0}{r}) = - \frac{1}{2\pi} \,  \frac{(12 x_0^2 -3 r^2) \ux}{\dkw^{7/2}}
$$
and for $k=4$
$$
B_{-4} =   \frac{1}{2\pi} \, 3! \,  \frac{3\cdot5\cdot7}{6\cdot7\cdot8} \, \frac{r^3 \ux}{\dkw^{9/2}} \, i^3 \, C_3^{-4}(i \frac{x_0}{r}) = \frac{1}{2\pi} \,  \frac{(60 x_0^3 - 45 x_0 r^2) \ux}{\dkw^{9/2}}
$$
Their distributional limits for $x_0 \rightarrow 0+$ are given by
$$
\left \{ \begin{array}{rcl}
a_{-2\ell} & = & (- \pux)^{2\ell-1} H  =   (-1)^{\ell-1} \, \frac{1}{2\pi} \, (2\ell-1)!! \, (2\ell-1)!! \, {\rm Fp} \displaystyle\frac{1}{r^{2\ell+1}}  \\[4mm]
b_{-2\ell} & = & (- \pux)^{2\ell-1} \delta       \end{array} \right .
$$
and
$$
\left \{ \begin{array}{rcl}
a_{-2\ell-1} & = &  \pux^{2\ell} \delta\\[2mm]
b_{-2\ell-1} & = & \pux^{2\ell} H
 = (-1)^{\ell-1} \, \frac{1}{2\pi} \, (2\ell-1)!! \, (2\ell+1)!! \,  {\rm Fp} \displaystyle\frac{1}{r^{2\ell+2}} \, \uom  \end{array} \right .
$$
They show the following properties, in fact valid regardless the dimension.

\begin{property}
\label{lem2}
One has for $j,k=1,2,\ldots$
\begin{itemize}
\item[(i)] $a_{-k} \xrightarrow{\hspace*{1mm} -\pux \hspace*{1mm}} b_{-k-1} \xrightarrow{\hspace*{1mm} -\pux \hspace*{1mm}} a_{-k-2}$
\item[(ii)] $\mathcal{H} \left [ a_{-k} \right ] = b_{-k}$, $\mathcal{H} \left [ b_{-k} \right ] = a_{-k}$
\item[(iii)] $a_{-j} \ast a_{-k} = a_{-j-k+1}$ \\
$a_{-j} \ast b_{-k} = b_{-j} \ast a_{-k} = b_{-j-k+1}$ \\
$b_{-j} \ast b_{-k} = a_{-j-k+1}$.
\end{itemize}
\end{property}

\subsection{The upstream potentials in $\mR_+^{3}$}

Let us have a look at the so--called {\em upstream} potentials. To start with the fundamental solution of the Laplace operator $\Delta_{3}$ in $\mR^{3}$, sometimes called Green's function, and here denoted by $\frac{1}{2}A_0(x_0,\ux)$, is given by
$$
\frac{1}{2}A_0(x_0,\ux) = - \frac{1}{4\pi} \frac{1}{\dkw^\onehalf}
$$
Its conjugate harmonic in $\mR^{3}_+$, in the sense of \cite{fbrdfs}, is
\begin{equation}
B_0(x_0,\ux) = \frac{2}{\sigma_{3}} \, \frac{\ux}{|\ux|^2} \, \mF_2 \left ( \frac{r}{x_0} \right )
\label{B0}
\end{equation}
where
$$
\mF_2(v) = \int_0^v \frac{\eta}{(1+\eta^2)^\frac{3}{2}} \, d\eta = \frac{v^2}{1+v^2+\sqrt{1+v^2}}
$$
leading to
$$
B_0(x_0,\ux) = \frac{1}{2\pi} \, \frac{\ux}{ \sqrt{x_0^2+r^2}\left(x_0+\sqrt{x_0^2+r^2}\right)}
$$
It is verified that
$$
\po B_0(x_0,\ux) = - \frac{1}{2\pi} \, \frac{\ux}{\dkw^{3/2}} = Q(x_0,\ux)
$$
and
$$
\pux  B_0(x_0,\ux) = - \frac{1}{2\pi} \, \frac{x_0}{\dkw^{3/2}} = - P(x_0,\ux)
$$
and also
$$
\lim_{x_0 \rightarrow 0+} \, B_0(x_0,\ux) = \frac{1}{2\pi} \, \frac{\ux}{r^2} = b_0(\ux)
$$
Note that $b_0(\ux) = \frac{1}{2\pi} \, \frac{\ux}{r^2} = - E_1$, with $E_1$ the fundamental solution of the Dirac operator $\pux^1 = \pux$.\\[2mm]
Green's function $A_0(x_0,\ux)$ itself shows the distributional limit
$$
\lim_{x_0 \rightarrow 0+} A_0(x_0,\ux) = - \frac{1}{2\pi} \, \frac{1}{r} = a_0(\ux)
$$
Note that $a_0(\ux) = - F_1$, $F_1 = \frac{1}{2\pi} \, \frac{1}{r}$ being the fundamental solution to the so--called Hilbert--Dirac operator $^1\mcH = (-\Delta_2)^\onehalf$ (see \cite{bdbds2, hidi}).\\

It is readily seen that $\overline{D} A_0 =  \overline{D} \overline{e_0} B_0 =  C_{-1}$. So $A_0(x_0,\ux)$ and $\overline{e_0} B_0(x_0,\ux)$ are conjugate harmonic potentials with respect to the operator $\overline{D}$, of the Cauchy kernel $C_{-1}(x_0,\ux)$ in $\mR^{3}_+$. Putting $C_0(x_0,\ux) = \frac{1}{2} A_0(x_0,\ux) + \frac{1}{2} \overline{e_0} B_0 (x_0,\ux)$, it follows that also $\overline{D} C_0 (x_0,\ux)  = C_{-1}(x_0,\ux)$,
which implies that $C_0(x_0,\ux)$ is a monogenic potential (or monogenic primitive) of the Cauchy kernel $C_{-1}(x_0,\ux)$ in $\mR^{3}_+$. 
Their distributional boundary values are intimately related, as mentioned in the following property, valid in general.
\begin{property}
\label{lemintiem}
One has
\begin{itemize}
\item[(i)] $-\pux a_0 = b_{-1} = H$; \quad $-\pux b_0 = a_{-1} = \delta$
\item[(ii)] $\mathcal{H} \left [a_0 \right ] = b_0$; \quad $\mathcal{H} \left [b_0 \right ] = a_0$
\end{itemize}
\end{property}

\begin{remark}
In the upper half of the complex plane the function $\ln(z)$ is a holomorphic potential (or primitive) of the Cauchy kernel $\frac{1}{z}$ and its real and imaginary components are the fundamental solution $\ln |z|$ of the Laplace operator, and its conjugate harmonic $i\, {\rm arg}(z)$ respectively. By similarity we could say that $C_0(x_0,\ux) = \frac{1}{2} A_0(x_0,\ux) + \frac{1}{2} \overline{e_0} B_0(x_0,\ux)$, being a monogenic potential of the Cauchy kernel $C_{-1}(x_0,\ux)$ and the sum of the fundamental solution $A_0(x_0,\ux)$ of the Laplace operator and its conjugate harmonic $\overline{e_0} B_0(x_0,\ux)$, is a {\em monogenic logarithmic function} in the upper half--space $\mR^{3}_+$.
\end{remark}

The construction of the sequence of {\em upstream} harmonic and monogenic potentials in $\mR^{3}_+$ is continued as follows.\\
The general expression for $A_1(x_0,\ux)$ established in \cite{bdbds1}, is not valid for $m=2$.  By direct calculation we obtain
$$
A_1(x_0,\ux) = - \frac{1}{2\pi} \, \ln{\left(x_0+\sqrt{\dk}\right)} 
$$
and it is verified that $\po A_1 = A_0$, $-\pux A_1 = B_0$ and $\lim_{x_0 \rightarrow 0+} \, A_1 = - \frac{1}{2\pi} \ln{r} = a_1(\ux)$.
Note that $a_1(\ux) = - \frac{1}{2\pi} \ln{r}  = E_2$ is the fundamental solution of the negative Laplace operator $\pux^2 = - \Delta_2$.\\

For its conjugate harmonic in $\mR^3_+$ we obtain
$$
B_1(x_0,\ux) = \frac{1}{2\pi} \, \frac{x_0 \ux}{r^2} \, \mF_2\left( \frac{r}{x_0}\right) -  \frac{1}{2\pi} \, \frac{\ux}{\dkw^\onehalf} =  \frac{1}{2\pi} \, \frac{\ux}{r^2} \left( x_0 - \sqrt{\dk}  \right)
$$
for which it is verified that $\po B_1 = B_0$, $-\pux B_1 = A_0$ and $\lim_{x_0 \rightarrow 0+} \, B_1 = - \frac{1}{2\pi} \frac{\ux}{r} = b_1(\ux)$.\\[2mm] 
Note that $b_1(\ux) =  - \frac{1}{2\pi} \frac{\ux}{r} = - \frac{1}{2\pi} \uom  = F_2$ is the fundamental solution of the operator $^2\mcH$ (see \cite{bdbds2}).\\

It follows that $\overline{D} A_{-1} =  \overline{D} \overline{e_0} B_{-1} =  C_{0}$, whence $A_1(x_0,\ux)$ and $B_1(x_0,\ux)$ are conjugate harmonic potentials in $\mR^{3}_+$ of the function $C_0(x_0,\ux)$ and
$$
C_1(x_0,\ux) = \frac{1}{2} A_1(x_0,\ux) + \frac{1}{2} \overline{e_0} B_1(x_0,\ux)
$$
is a monogenic potential in $\mR^{3}_+$ of $C_0$.
The above mentioned distributional boundary values show the following properties, valid in general.
\begin{property}
\label{lem54}
\rule{0mm}{0mm}
\begin{itemize}
\item[(i)] $- \pux a_1 = b_0$, $-\pux b_1 = a_0$
\item[(ii)] $\mathcal{H} \left [ a_1 \right ] = b_1$, $\mathcal{H} \left [ b_1 \right ] = a_1$
\end{itemize}
\end{property}

In the next step the general expressions for $A_2(x_0,\ux)$ and $B_2(x_0,\ux)$ are not valid. A direct computation yields
$$
A_2(x_0,\ux) = \frac{1}{2\pi} \, \left( \sqrt{\dk} - x_0 \ln{\left(x_0 + \sqrt{\dk}\right)}  \right)
$$
and it is verified that $-\pux A_2 = B_1$ and $\lim_{x_0 \rightarrow 0+} \, A_2 =  \frac{1}{2\pi} r = a_2(\ux)$.\\
Note that $a_2(\ux) = - F_3$, with $F_3 = - \frac{1}{2\pi} r$ the fundamental solution of the operator $^3\mcH$ (see \cite{bdbds2}).\\

For its conjugate harmonic in $\mR^3_+$ we find
$$
B_2(x_0,\ux) = \frac{\ux}{4\pi} \, \left(  \onehalf - \frac{x_0}{x_0 + \sqrt{\dk}} - \ln{\left( x_0 + \sqrt{\dk}  \right)}   \right)
$$
for which it is verified that $-\pux B_2 = A_1$ and $\lim_{x_0 \rightarrow 0+} \, B_2 =  \frac{\ux}{4\pi} \, (- \ln r + \onehalf) = b_2(\ux)$.
Note that $b_2(\ux) = - E_3$, with $E_3 = (\frac{1}{4\pi} \ln{r} - \frac{1}{8\pi})\ux$ the fundamental solution of the operator $\pux^3$.\\

It follows that
$$
C_2(x_0,\ux) = \frac{1}{2} A_2(x_0,\ux) + \frac{1}{2} \overline{e_0} B_2(x_0,\ux)
$$
is a monogenic potential in $\mR^{3}_+$ of $C_1$. The distributional limits show the following properties, also valid in general.
\begin{property}
\label{lem56}
\rule{0mm}{0mm}
\begin{itemize}
\item[(i)] $-\pux a_2 = b_1$, $- \pux b_2 = a_1$
\item[(ii)] $\mathcal{H} \left [ a_2 \right ] = b_2$, $\mathcal{H} \left [ b_2 \right ] = a_2$
\end{itemize}
\end{property}

Inspecting the above expressions for the harmonic potentials $A_1(x_0,\ux)$ and $A_2(x_0,\ux)$ we can put forward a general form for the potentials $A_j(x_0,\ux), j=1,2,\ldots$

\begin{proposition}
For $j=1,2,\ldots$ one has
$$
2\pi A_j(x_0,\ux) = P_j(x_0,r^2) \, \ln{(x_0 + \sqrt{x_0^2+r^2})} + Q_j(x_0,r^2) \,  \sqrt{x_0^2+r^2} + S_j(x_0,r^2)
$$
with
$$
P_{2k}(x_0,r^2) = p_{2k}^{2k-1} x_0^{2k-1} + p_{2k}^{2k-3} r^2 x_0^{2k-3} + \cdots + p_{2k}^{1} r^{2k-2}x_0
$$
$$
P_{2k+1}(x_0,r^2) = p_{2k+1}^{2k} x_0^{2k} + p_{2k+1}^{2k-2} r^2 x_0^{2k-2} + \cdots + p_{2k+1}^{0} r^{2k}
$$
and
$$
Q_{2k}(x_0,r^2) = q_{2k}^{2k-2} x_0^{2k-2} + q_{2k}^{2k-4} r^2 x_0^{2k-4} + \cdots + q_{2k}^{0} r^{2k-2}
$$
$$
Q_{2k+1}(x_0,r^2) = q_{2k+1}^{2k-1} x_0^{2k-1} + q_{2k+1}^{2k-3} r^2 x_0^{2k-3} + \cdots + q_{2k+1}^{1} r^{2k-2}x_0
$$
and
$$
S_{2k}(x_0,r^2) = s_{2k}^{2k-3} r^2 x_0^{2k-3} + s_{2k}^{2k-5} r^4 x_0^{2k-5} + \cdots + s_{2k}^{1} r^{2k-2} x_0
$$
$$
S_{2k+1}(x_0,r^2) = s_{2k+1}^{2k-2} r^2 x_0^{2k-2} + s_{2k+1}^{2k-4} r^4 x_0^{2k-4} + \cdots + s_{2k+1}^{0} r^{2k}
$$
all the coefficients $p_{2k}^j$, $p_{2k+1}^j$, $q_{2k}^j$, $q_{2k+1}^j$, $s_{2k}^j$ and $s_{2k+1}^j$ being real constants.
\end{proposition}

\pf
The harmonic potentials $A_1(x_0,\ux)$ and $A_2(x_0,\ux)$ computed above, fit into this general form. Now we will show that it is possible to determine unambiguously all the coefficients in the expression of $A_j(x_0,\ux)$ in terms of the coefficients in $A_{j-1}(x_0,\ux)$. To that end we impose on $A_j(x_0,\ux)$ the following two conditions, in line with its definition:
$$
\po (2\pi A_j(x_0,\ux)) = 2\pi A_{j-1}(x_0,\ux)
$$
and 
$$
\lim_{x_0 \rightarrow 0+} (2\pi A_j(x_0,\ux)) = 2\pi a_{j}(\ux)
$$
In the case where $j$ is even, say $j=2k$, this leads to the equations
\begin{equation}
\label{evenP}
\po P_{2k} = P_{2k-1}
\end{equation}
\begin{equation}
\label{evenS}
\po S_{2k} = S_{2k-1}
\end{equation}
\begin{equation}
\label{evenQ}
P_{2k} + x_0 Q_{2k} + (x_0^2+r^2) \, \po Q_{2k} = (x_0^2+r^2) \, \po Q_{2k-1}
\end{equation}
\begin{equation}
\label{evenlim}
q_{2k}^0 r^{2k-1} + s_{2k}^0 r^{2k-1} = 2\pi a_{2k}(\ux)
\end{equation}

From (\ref{evenP}) all the $p_{2k}$--coefficients may be determined as
$$
p_{2k}^j = \frac{1}{j} \, p_{2k-1}^{j-1}
$$
while from (\ref{evenS}) all the $s_{2k}$--coefficients follow by the similar relation
$$
s_{2k}^j = \frac{1}{j} \, s_{2k-1}^{j-1}
$$
Then all the $q_{2k}$--coefficients follow recursively from (\ref{evenQ}), and equation (\ref{evenlim}) may be used as a check.

In the case where $j$ is odd, say $j=2k+1$, the similar equations read
\begin{equation}
\label{oddP}
\po P_{2k+1} = P_{2k}
\end{equation}
\begin{equation}
\label{oddS}
\po S_{2k+1} = S_{2k}
\end{equation}
\begin{equation}
\label{oddQ}
P_{2k+1} + x_0 Q_{2k+1} + (x_0^2+r^2) \, \po Q_{2k+1} = (x_0^2+r^2) \, \po Q_{2k}
\end{equation}
and
\begin{equation}
\label{oddlim}
p_{2k+1}^0 r^{2k} \ln{r}+ s_{2k+1}^0 r^{2k} = 2\pi a_{2k+1}(\ux)
\end{equation}

From (\ref{oddP}) all the $p_{2k+1}$--coefficients, except $p_{2k+1}^0$, may be determined as
$$
p_{2k+1}^j = \frac{1}{j} \, p_{2k}^{j-1}
$$
while from (\ref{oddS}) all the $s_{2k+1}$--coefficients, except $s_{2k+1}^0$, follow by the similar relation
$$
s_{2k+1}^j = \frac{1}{j} \, s_{2k}^{j-1}
$$
The remaining coefficients $p_{2k+1}^0$ and $s_{2k+1}^0$ follow from (\ref{oddlim}).
All the $q_{2k+1}$--coefficients then follow recursively from (\ref{oddQ}).
\qed

\begin{remark}
It is obvious that solving equations (\ref{evenlim}) and (\ref{oddlim}) requires the knowledge of the distributional boundary values $a_j(\ux), j=1,2,\ldots$ There holds (see \cite{bdbds2})
$$
a_{2k}(\ux) = (-1)^{k+1} \frac{1}{2\pi} \frac{1}{((2k-1)!!)^2} \, r^{2k-1} \quad  (k = 1,2,\ldots)
$$
and
$$
a_{2k+1}(\ux) = (\alpha_{2k} \ln{r} + \beta_{2k}) \frac{\pi^{k+1}}{k!} \, r^{2k} \quad  (k = 0,1,2,\ldots)
$$
where the coefficients $\alpha_{2k}$ and $\beta_{2k}$ are defined recursively by
$$
\left \{ \begin{array}{rcl}
\alpha_{2k+2} & =  &  - \displaystyle\frac{1}{2\pi} \displaystyle\frac{1}{2k+2} \, \alpha_{2k}\\[4mm]
\beta_{2k+2} & = & - \displaystyle\frac{1}{2\pi} \displaystyle\frac{1}{2k+2} \, (\beta_{2k} - \displaystyle\frac{1}{k+1} \, \alpha_{2k})
\end{array} \right . \qquad k=0,1,2,\ldots 
$$
with starting values $\alpha_{0} =  - \displaystyle\frac{1}{2\pi^2}$ and $\beta_{0} = 0$, leading to their closed form
$$
\alpha_{2k} = \displaystyle\frac{(-1)^{k+1}}{2^{2k+1} \pi^{k+2} k!} \quad {\rm and } \quad 
\beta_{2k} =  \displaystyle\frac{(-1)^{k} H_{k}}{2^{2k+1} \pi^{k+2} k!}
$$
with $H_k = \sum_{n=1}^k \, \frac{1}{n}$.

\end{remark}

Now by Proposition 3.1 all the upstream $A_j$--potentials may be calculated recursively. Using the shorthands
$$
LOG = \ln{(x_0+\sqrt{x_0^2+r^2})} \quad {\rm and} \quad SQRT = \sqrt{x_0^2+r^2} 
$$
the outcome of our calculations is the following:
$$
2\pi A_3(x_0,\ux) = (-\onehalf x_0^2 + \frac{1}{4} r^2) LOG + \frac{3}{4} x_0 SQRT - \frac{1}{4} r^2
$$
$$
2\pi A_4(x_0,\ux) = (-\frac{1}{6} x_0^3 + \frac{1}{4} r^2 x_0) LOG + (\frac{11}{36} x_0^2 - \frac{1}{9} r^2) SQRT - \frac{1}{4} r^2 x_0
$$
$$
2\pi A_5(x_0,\ux) = (-\frac{1}{24} x_0^4 + \frac{1}{8} r^2 x_0^2 - \frac{1}{64} r^4) LOG + (\frac{25}{288} x_0^3 - \frac{55}{576} r^2 x_0) SQRT - \frac{1}{8} r^2 x_0^2 +  \frac{3}{128} r^4
$$
$$
2\pi A_6(x_0,\ux) = (-\frac{1}{120} x_0^5 + \frac{1}{24} r^2 x_0^3 - \frac{1}{64} r^4 x_0) LOG + (\frac{137}{7200} x_0^4 - \frac{607}{14400} r^2 x_0^2 + \frac{1}{225}  r^4) SQRT - \frac{1}{24} r^2 x_0^3 +  \frac{3}{128} r^4 x_0
$$

Note that once an upstream $A_j$--potential is determined, the $B_{j-1}$--potential follows readily by
$$
(-\pux)A_j(x_0,\ux) = B_{j-1}(x_0,\ux)
$$


\section{Harmonic and monogenic potentials in $\mR_+^{4}$}
\label{potentials4}


In this section we calculate the harmonic and monogenic potentials in upper half--space $\mR^{4}_+$. Quite naturally the structure of this section is completely similar to the foregoing one.

\subsection{The Cauchy kernel in $\mR_+^{4}$}

The starting point is the Cauchy kernel, i.e. the fundamental solution of the generalized Cauchy--Riemann operator $D$:
$$
C_{-1}(x_0,\ux) = \frac{1}{\sigma_{4}} \, \frac{x \overline{e_0}}{|x|^{4}} =  \frac{1}{2\pi^2} \, \frac{x_0 - \overline{e_0} \ux}{\dkwo^{2}}
$$
where we have put now $|\ux| = \rho$.
It may be decomposed in terms of the traditional Poisson kernels in $\mR^{4}_+$:
$$
C_{-1}(x_0,\ux) = \frac{1}{2} A_{-1}(x_0,\ux) + \frac{1}{2} \overline{e_0} \, B_{-1}(x_0,\ux)
$$
where, also mentioning the usual notations, for $x_0 >0$,
\begin{eqnarray*}
A_{-1}(x_0,\ux) & = & P(x_0,\ux) \ = \ \phantom{-} \frac{1}{\pi^2} \, \frac{x_0}{\dkwo^{2}} \label{A-1}\\
B_{-1}(x_0,\ux) & = & Q(x_0,\ux) \ = \ - \frac{1}{\pi^2} \, \frac{\ux}{\dkwo^{2}} \label{B-1}
\end{eqnarray*}
Their distributional limits for $x_0 \rightarrow 0+$ are given by
\begin{eqnarray*}
a_{-1}(\ux) & = & \lim_{x_0 \rightarrow 0+} A_{-1}(x_0,\ux) \ = \ \delta(\ux)\\
b_{-1}(\ux) & = & \lim_{x_0 \rightarrow 0+} B_{-1}(x_0,\ux) \ = \ H(\ux) \ = \ - \frac{1}{\pi^2} \, \mbox{Pv} \frac{\ux}{\rho^4}  \ = \ - \frac{1}{\pi^2} \, \mbox{Pv} \frac{\uom}{\rho^3}
\end{eqnarray*}
Note that, as in general, both distributional boundary values are linked by the Hilbert transform  $\mathcal{H}$ in $\mR^3$ with the above convolution kernel $H(\ux)$:
\begin{eqnarray*}
\mathcal{H} \left [ a_{-1} \right ] & = & \mathcal{H} \left [ \delta \right ] \ = \ H \ast \delta \ = \ H \ = \ b_{-1} \\
\mathcal{H} \left [ b_{-1} \right ] & = & \mathcal{H} \left [ H \right ] \ = \ H \ast H \ = \ \delta \ = \ a_{-1}
\end{eqnarray*}
Note also that $a_{-1}(\ux) = \delta(\ux) = E_0$ can be seen as the fundamental solution of the identity operator $\pux^0 = {\bf 1}$, and that $b_{-1}(\ux) =  - \frac{1}{\pi^2} \, \mbox{Pv} \frac{\uom}{\rho^3} = H(\ux) = F_0$ is the fundamental solution of the Hilbert operator $^0\mcH = \mcH$.

\subsection{The downstream potentials in $\mR_+^{4}$}

The first in the sequence of the {\em downstream} potentials is the function $C_{-2}$ defined by
$$
\overline{D} C_{-1} = C_{-2} = \frac{1}{2} A_{-2} + \frac{1}{2} \overline{e_0} B_{-2}
$$
Clearly it is monogenic  in $\mR^{4}_+$, since $DC_{-2} = D \overline{D} C_{-1} = \frac{1}{4} \Delta_{m+1} C_{-1} = 0$. 
The definition itself of $C_{-2}(x_0,\ux)$ implies that it shows the monogenic potential (or primitive) $C_{-1}(x_0,\ux)$ and the conjugate harmonic potentials $A_{-1}(x_0,\ux)$ and $\overline{e_0} B_{-1}(x_0,\ux)$.\\
The harmonic component $A_{-2}(x_0,\ux)$ may be calculated as $\po A_{-1}$ or as $-\pux B_{-1}$, and even the general expression is valid, leading to
$$
A_{-2}(x_0,\ux) = \frac{1}{\pi^2} \, \frac{-3x_0^2 + \rho^2}{\dkwo^{3}}
$$
Similarly the harmonic component $B_{-2}(x_0,\ux)$ may be calculated as $\po B_{-1}$ or as $-\pux A_{-1}$, and also the general expression is valid, leading to
$$
B_{-2}(x_0,\ux) = \frac{4}{\pi^2} \, \frac{ x_0 \ux}{\dkwo^{3}}
$$
The distributional limits for $x_0 \rightarrow 0+$ of these harmonic potentials are given by
$$
\left \{ \begin{array}{rcl}
a_{-2}(\ux) = \lim_{x_0 \rightarrow 0+} A_{-2}(x_0,\ux) & = & \displaystyle\frac{1}{\pi^2} \, {\rm Fp} \displaystyle\frac{1}{\rho^4}\\[4mm]
b_{-2}(\ux) = \lim_{x_0 \rightarrow 0+} B_{-2}(x_0,\ux) & = & - \pux \delta
\end{array} \right .
$$
Note that $a_{-2}(\ux) = - F_{-1}$, with $F_{-1} = \pux H$ the fundamental solution of the operator $^{-1}\mcH$, while 
$b_{-2}(\ux) = - E_{-1}$, with $E_{-1} = \pux \delta$ the fundamental solution of the operator $\pux^{-1}$ (see \cite{bdbds2}).\\

The sequence of {\em downstream} monogenic potentials in $\mR_+^{m+1}$ is defined, as in general dimension, by
$$
C_{-k-1} = \overline{D} C_{-k} = \overline{D}^2 C_{-k+1} = \ldots = \overline{D}^k C_{-1}, \qquad k=1,2,\ldots
$$
where each monogenic potential decomposes into two conjugate harmonic potentials:
$$
C_{-k-1} = \frac{1}{2} A_{-k-1} + \frac{1}{2} \overline{e_0} B_{-k-1}, \qquad k=1,2,\ldots
$$
The harmonic components $A_{-k}(x_0,\ux)$ are real--valued and given by
$$
A_{-k}(x_0,\ux) = (-1)^{k+1} \, \frac{1}{\pi^2} \,  \frac{2^{k-1} (k!)^2}{(k+2)(k+3)\cdots(2k+1)} \, \frac{\rho^k}{\dkwo^{k+1}} \, i^k \, C_k^{-k-\onehalf}(i \frac{x_0}{\rho})
$$
We have explicitly calculated $A_{-k}$ for the $k$--values $1,2,3$ and $4$. We obtain\\

\noindent
for $k=1$
$$
A_{-1} = \frac{1}{\pi^2} \, \frac{1}{3} \, \frac{\rho}{\dkwo^{2}} \, i \, C_1^{-3/2}(i \frac{x_0}{\rho}) = \frac{1}{\pi^2} \,  \frac{x_0}{\dkwo^{2}}
$$
for $k=2$
$$
A_{-2} = - \frac{2}{\pi^2} \, \frac{2\cdot2}{4\cdot5} \, \frac{\rho^2}{\dkwo^{3}} \, i^2 \, C_2^{-5/2}(i \frac{x_0}{\rho}) = \frac{1}{\pi^2} \,  \frac{-3 x_0^2 + \rho^2}{\dkwo^{3}}
$$
for $k=3$
$$
A_{-3} =  \frac{4}{\pi^2} \, \frac{6\cdot6}{5\cdot6\cdot7} \, \frac{\rho^3}{\dkwo^{4}} \, i^3 \, C_3^{-7/2}(i \frac{x_0}{\rho}) = \frac{1}{\pi^2} \,  \frac{12 x_0^3 - 12 x_0 \rho^2}{\dkwo^{4}}
$$
and for $k=4$
$$
A_{-4} = - \frac{8}{\pi^2} \, \frac{24\cdot24}{6\cdot7\cdot8\cdot9} \, \frac{\rho^4}{\dkwo^{5}} \, i^4 \, C_4^{-9/2}(i \frac{x_0}{\rho}) = \frac{1}{\pi^2} \,  \frac{-60 x_0^4 + 120  x_0^2  \rho^2 - 12 \rho^4}{\dkwo^{5}}
$$
The conjugate harmonic components $B_{-k}(x_0,\ux)$ are Clifford vector--valued and given by
$$
B_{-k}(x_0,\ux) = (-1)^{k} \, \frac{1}{\pi^2} \, \frac{2^{k-1} (k-1)! k!}{(k+3)(k+4)\cdots(2k+1)} \, \frac{\rho^{k-1} \ux}{\dkwo^{k+1}} \, i^{k-1} \, C_{k-1}^{-k-\onehalf}(i \frac{x_0}{\rho})
$$
We have explicitly calculated $B_{-k}$ for the $k$--values $1,2,3$ and $4$. We obtain\\

\noindent
for $k=1$
$$
B_{-1} = - \frac{1}{\pi^2} \, \frac{\ux}{\dkwo^{2}} \, C_0^{-3/2}(i \frac{x_0}{\rho}) = - \frac{1}{\pi^2} \,  \frac{\ux}{\dkwo^{2}}
$$
for $k=2$
$$
B_{-2} =  \frac{2}{\pi^2} \, \frac{1! 2!}{5} \, \frac{\rho \ux}{\dkw^{3}} \, i \, C_1^{-5/2}(i \frac{x_0}{\rho}) = \frac{4}{\pi^2} \,  \frac{x_0 \ux}{\dkwo^{3}}
$$
for $k=3$
$$
B_{-3} =  - \frac{4}{\pi^2} \,  \frac{2!3!}{6\cdot7} \, \frac{\rho^2 \ux}{\dkw^{4}} \, i^2 \, C_2^{-7/2}(i \frac{x_0}{\rho}) =  \frac{4}{\pi^2} \,  \frac{(-5 x_0^2 + \rho^2) \ux}{\dkwo^{4}}
$$
and for $k=4$
$$
B_{-4} =   \frac{8}{\pi^2}  \,  \frac{3!4!}{7\cdot8\cdot9} \, \frac{\rho^3 \ux}{\dkw^{5}} \, i^3 \, C_3^{-9/2}(i \frac{x_0}{\rho}) = \frac{24}{\pi^2} \,  \frac{(5 x_0^3 - 3 x_0 \rho^2) \ux}{\dkwo^{5}}
$$
Their distributional limits for $x_0 \rightarrow 0+$ are given by
$$
\left \{ \begin{array}{rcl}
a_{-2\ell} & = & (- \pux)^{2\ell-1} H  =   (-1)^{\ell-1}  2^{\ell-1} (2\ell-1)!! \, \ell! \, \displaystyle\frac{1} {\pi^2} \, {\rm Fp} \displaystyle\frac{1}{\rho^{2\ell+2}}  \\[4mm]
b_{-2\ell} & = & (- \pux)^{2\ell-1} \delta       \end{array} \right .
$$
and
$$
\left \{ \begin{array}{rcl}
a_{-2\ell-1} & = &  \pux^{2\ell} \delta\\[2mm]
b_{-2\ell-1} & = & \pux^{2\ell} H
 = (-1)^{\ell-1}  2^{\ell} (2\ell-1)!! \, (\ell+1)! \, \displaystyle\frac{1}{\pi^2} \,  {\rm Fp} \displaystyle\frac{1}{\rho^{2\ell+3}} \, \uom  \end{array} \right .
$$
They show the by now traditional properties (see Property 2).

\subsection{The upstream potentials in $\mR_+^{4}$}

For the so--called {\em upstream} potentials, we start with the fundamental solution of the Laplace operator $\Delta_{4}$ in $\mR^{4}$, denoted by $\frac{1}{2}A_0(x_0,\ux)$, and given by
$$
\frac{1}{2}A_0(x_0,\ux) = - \frac{1}{4\pi^2} \frac{1}{\dko}
$$
Its conjugate harmonic in $\mR^{4}_+$, in the sense of \cite{fbrdfs}, is
\begin{equation}
B_0(x_0,\ux) = \frac{2}{\sigma_{4}} \, \frac{\ux}{\rho^3} \, \mF_3 \left ( \frac{\rho}{x_0} \right )
\label{B0}
\end{equation}
where
$$
\mF_3(v) = \int_0^v \frac{\eta^2}{(1+\eta^2)^2} \, d\eta = \onehalf \left( \arctan{v} - \frac{v}{1+v^2}\right)
$$
leading to
$$
B_0(x_0,\ux) = \frac{1}{2\pi^2} \, \frac{\ux}{\rho^3} \, \left( \arctan{\frac{\rho}{x_0}}  - \frac{x_0 \rho}{\dko}  \right)
$$
It is verified that
$$
\po B_0(x_0,\ux) = - \frac{1}{\pi^2} \, \frac{\ux}{\dkwo^{2}} = Q(x_0,\ux)
$$
and
$$
\pux  B_0(x_0,\ux) = - \frac{1}{\pi^2} \, \frac{x_0}{\dko} = - P(x_0,\ux)
$$
and also
$$
\lim_{x_0 \rightarrow 0+} \, B_0(x_0,\ux) = \frac{1}{4\pi} \, \frac{\ux}{\rho^3} = b_0(\ux)
$$
Note that $b_0(\ux) = \frac{1}{4\pi} \, \frac{\ux}{\rho^3} = \frac{1}{4\pi} \, \frac{\uom}{\rho^2} = -E_1$, with $E_1$ the fundamental solution of the Dirac operator $\pux^1 = \pux$.\\[2mm]
Green's function $A_0(x_0,\ux)$ itself shows the distributional limit
$$
 \lim_{x_0 \rightarrow 0+} A_0(x_0,\ux) = - \frac{1}{2\pi^2} \, \frac{1}{\rho^2} = a_0(\ux)
$$
Note that $a_0(\ux) = - F_1$, $F_1 = \frac{1}{2\pi^2} \, \frac{1}{\rho^2}$ being the fundamental solution to the so--called Hilbert--Dirac operator $^1\mcH = (-\Delta_3)^\onehalf$ (see \cite{bdbds2, hidi}).\\

It is readily seen that $\overline{D} A_0 =  \overline{D} \overline{e_0} B_0 =  C_{-1}$. So $A_0(x_0,\ux)$ and $\overline{e_0} B_0(x_0,\ux)$ are conjugate harmonic potentials (or primitives), with respect to the operator $\overline{D}$, of the Cauchy kernel $C_{-1}(x_0,\ux)$ in $\mR^{3}_+$. Putting $C_0(x_0,\ux) = \frac{1}{2} A_0(x_0,\ux) + \frac{1}{2} \overline{e_0} B_0 (x_0,\ux)$, it follows that also $\overline{D} C_0 (x_0,\ux)  = C_{-1}(x_0,\ux)$,
which implies that $C_0(x_0,\ux)$ is a monogenic potential (or primitive) of the Cauchy kernel $C_{-1}(x_0,\ux)$ in $\mR^{3}_+$. 
Their distributional boundary values are intimately related, as shown in the similar Property 3.

\begin{remark}
As in any dimension we could again say that $C_0(x_0,\ux) = \frac{1}{2} A_0(x_0,\ux) + \frac{1}{2} \overline{e_0} B_0(x_0,\ux)$, being a monogenic potential of the Cauchy kernel $C_{-1}(x_0,\ux)$ and the sum of the fundamental solution $A_0(x_0,\ux)$ of the Laplace operator and its conjugate harmonic $\overline{e_0} B_0(x_0,\ux)$, is a {\em monogenic logarithmic function} in the upper half--space $\mR^{4}_+$.
\end{remark}

The construction of the sequence of {\em upstream} harmonic and monogenic potentials in $\mR^{4}_+$ is continued as follows.\\
The general expression for $A_1(x_0,\ux)$, established in \cite{bdbds1}, remains valid for $m=3$:
$$
A_1(x_0,\ux) =  \frac{1}{2\pi^2} \, \frac{1}{\rho} \, \arctan{\frac{\rho}{x_0}}
$$
and it is verified that $\po A_1 = A_0$, $-\pux A_1 = B_0$ and $\lim_{x_0 \rightarrow 0+} \, A_1 = \frac{1}{4\pi} \, \frac{1}{\rho} = a_1(\ux)$. Note that this distributional boundary value $a_1(\ux) = E_2$ is the fundamental solution of the  negative Laplace operator $\pux^2 = - \Delta_3$.\\
For its conjugate harmonic in $\mR^4_+$ we obtain by direct calculation
$$
B_1(x_0,\ux) = \frac{1}{2\pi^2} \, \frac{\ux}{\rho^2} \, \left( \frac{x_0}{\rho} \arctan{\frac{\rho}{x_0}} - 1 \right)
$$
for which it is verified that $\po B_1 = B_0$, $-\pux B_1 = A_0$ and $\lim_{x_0 \rightarrow 0+} \, B_1 = - \frac{1}{2\pi^2} \frac{\ux}{\rho^2} = b_1(\ux)$. Note that $b_1(\ux) = - \frac{1}{2\pi^2} \frac{\uom}{\rho} = F_2$ is the fundamental solution of the operator $^2\mcH$ (see \cite{bdbds2}).\\[2mm]

It follows that $\overline{D} A_{-1} =  \overline{D} \overline{e_0} B_{-1} =  C_{0}$, whence $A_1(x_0,\ux)$ and $B_1(x_0,\ux)$ are conjugate harmonic potentials in $\mR^{3}_+$ of the function $C_0(x_0,\ux)$ and
$$
C_1(x_0,\ux) = \frac{1}{2} A_1(x_0,\ux) + \frac{1}{2} \overline{e_0} B_1(x_0,\ux)
$$
is a monogenic potential in $\mR^{3}_+$ of $C_0$.
The above mentioned distributional boundary values show properties similar to those of Property 4.\\

In the next step the general expressions for $A_2(x_0,\ux)$ is not valid. A direct computation yields
$$
A_2(x_0,\ux) = \frac{1}{2\pi^2} \, \left( \frac{x_0}{\rho} \arctan{\frac{\rho}{x_0}} + \ln{ \sqrt{\dko}} \right)
$$
and it is verified that $-\pux A_2 = B_1$ and $\lim_{x_0 \rightarrow 0+} \, A_2 =  \frac{1}{2\pi^2} \ln{\rho} = a_2(\ux)$. 
Note that $a_2(\ux) = - F_3$, $F_3 = - \frac{1}{2\pi^2} \ln{\rho}$ being the fundamental solution of the operator $^3\mcH$ (see \cite{bdbds2}).

For its conjugate harmonic in $\mR^4_+$ we find
$$
B_2(x_0,\ux) = \frac{1}{4\pi^2} \, \left( \frac{\ux}{\rho^3} \dkwo \arctan{\frac{\rho}{x_0}}  -  \frac{\ux x_0}{\rho^2}  \right)
$$
for which it is verified that $\po B_2 = B_1$, $-\pux B_2 = A_1$ and $\lim_{x_0 \rightarrow 0+} \, B_2 =  \frac{1}{8\pi} \, \frac{\ux}{\rho} = b_2(\ux)$. Note that $b_2(\ux) = -E_3$, with $E_3$ the fundamental solution of $\pux^3$.\\
It follows that
$$
C_2(x_0,\ux) = \frac{1}{2} A_2(x_0,\ux) + \frac{1}{2} \overline{e_0} B_2(x_0,\ux)
$$
is a monogenic potential in $\mR^{3}_+$ of $C_1$. The distributional limits show the properties similar to those of Property 5.\\

Inspecting the above expressions for the harmonic potentials $A_1(x_0,\ux)$ and $A_2(x_0,\ux)$ we can put forward a general form for the potentials $A_j(x_0,\ux), j=1,2,\ldots$

\begin{proposition}
For $j=1,2,\ldots$ one has
$$
2\pi A_j(x_0,\ux) = U_j(x_0,r^2) \, \frac{1}{r} \arctan{\frac{r}{x_0}} + V_j(x_0,r^2) \,  \ln{\sqrt{x_0^2+r^2}} + W_j(x_0,r^2)
$$
with
$$
U_{2k}(x_0,r^2) = u_{2k}^{2k-1} x_0^{2k-1} + u_{2k}^{2k-3} r^2 x_0^{2k-3} + \cdots + u_{2k}^{1} r^{2k-2}x_0
$$
$$
U_{2k+1}(x_0,r^2) = u_{2k+1}^{2k} x_0^{2k} + u_{2k+1}^{2k-2} r^2 x_0^{2k-2} + \cdots + u_{2k+1}^{0} r^{2k}
$$
and
$$
V_{2k}(x_0,r^2) = v_{2k}^{2k-2} x_0^{2k-2} + v_{2k}^{2k-4} r^2 x_0^{2k-4} + \cdots + v_{2k}^{0} r^{2k-2}
$$
$$
V_{2k+1}(x_0,r^2) = v_{2k+1}^{2k-1} x_0^{2k-1} + v_{2k+1}^{2k-3} r^2 x_0^{2k-3} + \cdots + v_{2k+1}^{1} r^{2k-2}x_0
$$
and
$$
W_{2k}(x_0,r^2) = w_{2k}^{2k-2} x_0^{2k-2} + w_{2k}^{2k-4} r^2 x_0^{2k-4} + \cdots + w_{2k}^{0} r^{2k-2}
$$
$$
W_{2k+1}(x_0,r^2) = w_{2k+1}^{2k-1} x_0^{2k-1} + w_{2k+1}^{2k-3} r^2 x_0^{2k-3} + \cdots + w_{2k+1}^{1} r^{2k-2} x_0
$$
all the coefficients $u_{2k}^j$, $u_{2k+1}^j$, $v_{2k}^j$, $v_{2k+1}^j$, $w_{2k}^j$ and $w_{2k+1}^j$ being real constants.
\end{proposition}

\pf
The proof is similar to that of Proposition 3.1. The harmonic potentials $A_1(x_0,\ux)$ and $A_2(x_0,\ux)$ computed above, fit into this general form. Now we will show that it is possible to determine unambiguously all the coefficients in the expression of $A_j(x_0,\ux)$ in terms of the coefficients in $A_{j-1}(x_0,\ux)$. To that end we impose on $A_j(x_0,\ux)$ the following two conditions, in line with its definition:
$$
\po (2\pi A_j(x_0,\ux)) = 2\pi A_{j-1}(x_0,\ux)
$$
and 
$$
\lim_{x_0 \rightarrow 0+} (2\pi A_j(x_0,\ux)) = 2\pi a_{j}(\ux)
$$
In the case where $j$ is even, say $j=2k$, this leads to the equations
\begin{equation}
\label{evenU}
\po U_{2k} = U_{2k-1}
\end{equation}
\begin{equation}
\label{evenV}
\po V_{2k} = V_{2k-1}
\end{equation}
\begin{equation}
\label{evenW}
- U_{2k} + x_0 V_{2k} + (x_0^2+r^2) \, \po W_{2k} = (x_0^2+r^2) \, W_{2k-1}
\end{equation}
\begin{equation}
\label{evenlim3}
v_{2k}^0 r^{2k-2} \ln{r} + w_{2k}^0 r^{2k-2} = 2\pi a_{2k}(\ux)
\end{equation}

From (\ref{evenU}) all the $u_{2k}$--coefficients may be determined as
$$
u_{2k}^j = \frac{1}{j} \, u_{2k-1}^{j-1}
$$
while from (\ref{evenV}) all the $v_{2k}$--coefficients, except $v_{2k}^0$, follow by the similar relation
$$
v_{2k}^j = \frac{1}{j} \, v_{2k-1}^{j-1}
$$
The coefficients $v_{2k}^0$ and $w_{2k}^0$ follow from (\ref{evenlim3}) and
all the other $w_{2k}$--coefficients follow recursively from the system of equations (\ref{evenW}). The last equation in this system can be used to check the value of $v_{2k}^0$.\\

In the case where $j$ is odd, say $j=2k+1$, the similar equations read
\begin{equation}
\label{oddU}
\po U_{2k+1} = U_{2k}
\end{equation}
\begin{equation}
\label{oddV}
\po V_{2k+1} = V_{2k}
\end{equation}
\begin{equation}
\label{oddW}
- U_{2k+1} + x_0 V_{2k+1} + (x_0^2+r^2) \, \po W_{2k+1} = (x_0^2+r^2) \, W_{2k}
\end{equation}
and
\begin{equation}
\label{oddlim3}
\frac{\pi}{2} u_{2k+1}^0 r^{2k-1}  = 2\pi a_{2k+1}(\ux)
\end{equation}

From (\ref{oddU}) all the $u_{2k+1}$--coefficients, except $u_{2k+1}^0$, may be determined as
$$
u_{2k+1}^j = \frac{1}{j} \, u_{2k}^{j-1}
$$
The coefficient $u_{2k+1}^0$ follows from (\ref{oddlim3}).
From (\ref{oddV}) all the $v_{2k}$--coefficients follow by the similar relation
$$
v_{2k+1}^j = \frac{1}{j} \, v_{2k}^{j-1}
$$
All the $w_{2k+1}$--coefficients follow recursively from the system of equations (\ref{oddW}), and the last equation in this system can be used to check the value of $u_{2k+1}^0$.
\qed

\begin{remark}
It is obvious that solving equations (\ref{evenlim3}) and (\ref{oddlim3}) requires the knowledge of the distributional boundary values $a_j(\ux), j=1,2,\ldots$ There holds (see \cite{bdbds2})
$$
a_{2k+1}(\ux) = (-1)^{k} \frac{1}{\pi} \frac{1}{2^{k+2}} \frac{1}{k!} \frac{1}{(2k-1)!!} \, r^{2k-1} \quad  (k = 0,1,2,\ldots)
$$
and
$$
a_{2k}(\ux) = - (\alpha_{2k-2} \ln{r} + \beta_{2k-2}) \frac{2^k \pi^{k}}{(2k-1)!!} \, r^{2k-2} \quad  (k = 1,2,\ldots)
$$
where the coefficients $\alpha_{2k}$ and $\beta_{2k}$ are defined recursively by
$$
\left \{ \begin{array}{rcl}
\alpha_{2k+2} & =  &  - \displaystyle\frac{1}{2\pi} \displaystyle\frac{1}{2k+2} \, \alpha_{2k}\\[4mm]
\beta_{2k+2} & = & - \displaystyle\frac{1}{2\pi} \displaystyle\frac{1}{2k+2} \, (\beta_{2k} - \displaystyle\frac{4k+5}{(2k+2)(2k+3)} \, \alpha_{2k})
\end{array} \right . \qquad k=0,1,2,\ldots 
$$
with starting values $\alpha_{0} =  - \displaystyle\frac{1}{4\pi^3}$ and $\beta_{0} = 0$, leading to their closed form
$$
\alpha_{2k} = \displaystyle\frac{(-1)^{k+1}}{2^{2k+2}\pi^{k+3}k!} \quad {\rm and} \quad
\beta_{2k} =  \displaystyle\frac{(-1)^{k} (H_{2k+1}-1)}{2^{2k+2}\pi^{k+3}k!}
$$
\end{remark}

Now by Proposition 4.1 all the upstream $A_j$--potentials may be calculated recursively. Using the shorthands
$$
QUAT =  \frac{1}{r} \arctan{\frac{r}{x_0}}  \quad {\rm and} \quad LNSQ = \ln{\sqrt{x_0^2+r^2}} 
$$
the outcome of our calculations is the following:
$$
2\pi A_3(x_0,\ux) = \frac{1}{2\pi} (x_0^2 - r^2) QUAT + \frac{1}{\pi} x_0 LNSQ - \frac{1}{2\pi} x_0
$$
$$
2\pi A_4(x_0,\ux) = (\frac{1}{6\pi} x_0^3 - \frac{1}{2\pi} r^2 x_0) QUAT + (\frac{1}{2\pi} x_0^2 - \frac{1}{6\pi} r^2) LNSQ - \frac{5}{12\pi} x_0^2 + \frac{5}{36\pi} r^2
$$
$$
2\pi A_5(x_0,\ux) = (\frac{1}{24\pi} x_0^4 - \frac{1}{4\pi} r^2 x_0^2 + \frac{1}{24\pi} r^4) QUAT + (\frac{1}{6\pi} x_0^3 - \frac{1}{6\pi} r^2 x_0) LNSQ - \frac{13}{72\pi} x_0^3 +  \frac{13}{72\pi} r^2 x_0
$$
$$
2\pi A_6(x_0,\ux) = (\frac{1}{120\pi} x_0^5 - \frac{1}{12\pi} r^2 x_0^3 + \frac{1}{24\pi} r^4 x_0) QUAT + (\frac{1}{24\pi} x_0^4 - \frac{1}{12\pi} r^2 x_0^2 + \frac{1}{120\pi}  r^4) LNSQ
$$
$$
 - \frac{77}{1440\pi} x_0^4 +  \frac{77}{720\pi} r^2 x_0^2 - \frac{77}{7200\pi} r^4
$$

\vspace*{3mm}
Note that once an upstream $A_j$--potential is determined, the $B_{j-1}$--potential follows readily by
$$
(-\pux)A_j(x_0,\ux) = B_{j-1}(x_0,\ux)
$$


\section{Conclusion}

The problem of generalizing, within the framework of Clifford analysis, the holomorphic function $\ln{z}$ in the upper half of the complex plane to a monogenic logarithmic function in the upper half of Euclidean space $\mR{m+1}_+$, led to a doubly infinite chain of monogenic and conjugate harmonic potentials or primitives with the Cauchy--kernel as a pivot element (see \cite{bdbds1}. Their distributional boundary limits in $\mR^m$ turned out to be the fundamental solutions of positive and negative integer powers of the Dirac operator $\pux$ and the Hilbert--Dirac operator, i. e. a convolution operator with kernel $\pux H$, $H$ being the multidimensional Hilbert kernel (see \cite{hidi}). These operators are special cases of the operators $^\mu \pux$ and $^\mu H$ defined for the complex parameter $\mu \in \mC$ (see \cite{disturb, bdbds2}). Those results depend, quite naturally, on the dimension $m$ of the Euclidean space considered. There is the traditional phenomenon in Clifford analysis that the parity of the dimension has a substantial impact, one could even speak of an {\em even} and and {\em odd} world with diverging results. Next to that there is the problem that for specific dimensions $m=2$ and $m=3$, which are in fact the most important dimensions with a view on application, the general formulae are no longer valid. Both cases need fresh ad hoc calculations, and the obtained low dimensional results are reported on in the underlying paper.




\begin{thebibliography}{06}

\bibitem{bdbds1} F.\ Brackx, H.\ De Bie, H.\ De Schepper, On a Chain of Harmonic and Monogenic Potentials in Euclidean Half--space, arXiv:1210.2044 (2012), submitted.

\bibitem{bdbds2} F.\ Brackx, H.\ De Bie, H.\ De Schepper, Distributional Boundary Values of Harmonic Potentials in Euclidean half--space as Fundamental Solutions of Convolution Operators in Clifford Analysis. In: G.\ Gentili, I.\ Sabadini, M.\ Shapiro, F.\ Sommen and D.C.\ Struppa, \textit{Advances in Hypercomplex Analysis}, Springer INdAM Series 1, Springer Verlag (2012), 15--37.

\bibitem{distrib} F.\ Brackx, B.\ De Knock, H.\ De Schepper, D.\ Eelbode, A Calculus Scheme for Clifford Distributions, \textit{Tokyo J. Math.} {\bf 29(2)} (2006), 495--513.

\bibitem{red} F.\ Brackx, R.\ Delanghe, F.\ Sommen, \textit{Clifford Analysis}, Pitman Publishers (Boston--London--Melbourne, 1982).

\bibitem{fbrdfs} F.\ Brackx, R.\ Delanghe, F.\ Sommen, On Conjugate Harmonic Functions in Euclidean Space, \textit{Math. Meth. Appl. Sci.} {\bf 25} (2002), 1553--1562.

\bibitem{fb1} F.\ Brackx, R.\ Delanghe, F.\ Sommen, Spherical means and distributions in Clifford analysis.
In: T.\ Qian, Th.\ Hempfling, A.\ McIntosh, F.\ Sommen (eds), \textit{Advances in Analysis and Geometry: New Developments Using Clifford Algebra}, \textit{Trends in Mathematics}, Birkh\"auser (Basel, 2004).

\bibitem{fb2} F.\ Brackx, R.\ Delanghe, F.\ Sommen, Spherical means, distributions and convolution operators in Clifford analysis, {\em Chin. Ann. Math.} {\bf 24B(2)} (2003), 133--146.

\bibitem{hidi}	F.\ Brackx, H.\ De Schepper, Hilbert-Dirac Operators in Clifford Analysis, \textit{Chin. Ann. Math.} {\bf 26B(1)} (2005), 1--14.

\bibitem{green} R.\ Delanghe, F.\ Sommen, V.\ Sou\v{c}ek, \textit{Clifford algebra and spinor-valued functions -- A function theory for the Dirac operator}, Kluwer Academic Publishers (Dordrecht, 1992).

\bibitem{gilmur}J.\ Gilbert, M.\ Murray, \textit{Clifford Algebra and Dirac Operators in Harmonic Analysis}, Cambridge University Press (Cambridge, 1991).

\bibitem{slang} S.\ Lang, \textit{Complex Analysis}, Graduate Texts in Mathematics, 103, Springer--Verlag (New York, 1999).

\bibitem{porteous} I.\ Porteous, \textit{Topological Geometry}, Van Nostrand Reinhold Company (London--New York--Toronto--Melbourne, 1969).


\end{thebibliography}
\end{document}